\documentclass{aims}
\usepackage{amsmath}
  \usepackage{paralist}
  \usepackage{graphics} 
  \usepackage{epsfig} 
\usepackage{graphicx}  \usepackage{epstopdf}
 \usepackage[colorlinks=true]{hyperref}

 \usepackage{amssymb}
\usepackage{stmaryrd}
\usepackage{amsfonts}
\usepackage{latexsym,amssymb,amsfonts,amsbsy,amsthm}
\usepackage[usenames]{color}
\usepackage{xspace,colortbl}
\hypersetup{urlcolor=blue, citecolor=red}

\def\currentvolume{X}
 \def\currentissue{X}
  \def\currentyear{200X}
   \def\currentmonth{XX}
    \def\ppages{X--XX}
     \def\DOI{10.3934/xx.xx.xx.xx}

\newtheorem{theorem}{Theorem}[section]

\newtheorem{lemma}[theorem]{Lemma}

\theoremstyle{definition}
\newtheorem{definition}[theorem]{Definition}

\title[Two coupled Maxwell systems in the temporal gauge]
      {Global solutions of two coupled Maxwell systems in the temporal gauge}

\author[Jianjun Yuan]{}


 \email{yuanjjn2@163.com}

\begin{document}
\maketitle

\centerline{\scshape Jianjun Yuan}
\medskip
{\footnotesize
 \centerline{The College of Information and Technology}
   \centerline{Nanjing University of Chinese Medicine}
   \centerline{Nanjing 210046, China}
} 

\medskip


\bigskip


\begin{abstract}
\noindent In this paper, we consider the Maxwell-Klein-Gordon and Maxwell-Chern-Simons-Higgs systems
in the temporal gauge. By using the fact that when the spatial gauge potentials are in the Coulomb gauge, their
$\dot{H}^1$ norms can be controlled by the energy of the corresponding system and their $L^2$ norms, and the gauge
invariance of the systems, we show that finite energy solutions of these two systems exist globally in this gauge.
\end{abstract}

\section{Introduction}
The Lagrangian density of the (3+1)-dimensional Maxwell-Klein-Gordon system
and the (2+1)-dimensional Maxwell-Chern-Simons-Higgs system are given respectively by
\begin{eqnarray}
\mathcal{L}_{MKG}=-\frac{1}{4}F^{\alpha\beta}F_{\alpha\beta}-\frac{1}{2}
D^{\mu}\phi\overline{D_{\mu}\phi},
\end{eqnarray}
and
\begin{eqnarray}
\mathcal{L}_{MCSH}=&-&\frac{1}{4}F^{\mu\nu}F_{\mu\nu}+\frac{\kappa}{4}\epsilon^{\mu\nu\rho}F_{\mu\nu}A_{\rho}+
D_{\mu}\phi\overline{D^{\mu}\phi}\nonumber\\
&+&\frac{1}{2}\partial_{\mu}N\partial^{\mu}N-\frac{1}{2}(e|\phi|^2+
\kappa N-ev^2)^2-e^2N^2|\phi|^2,
\end{eqnarray}
where $A_{\alpha}\in R$ is the gauge fields,  $\phi$ is a complex scalar
field, $N$ is a real scalar field, $F_{\alpha\beta}=\partial_{\alpha}A_{\beta}-\partial
_{\beta}A_{\alpha}$ is the curvature, $D_{\mu}=\partial_{\mu}-iA_{\mu}$ is the covariant derivative
for the MKG system and
$D_{\mu}=\partial_{\mu}-ieA_{\mu}$ is the covariant derivative for the MCSH system,
$e$ is the charge of the electron, $\kappa>0$ is the Chern-Simons
constant, $v$ is a nonzero constant, $\epsilon^{\mu\nu\rho}$ is the totally skew-symmetric tensor with $\epsilon^{012}=1$.
For the MKG system, indices are raised and lowered with respect to the Minkowski metric $g_{\alpha\beta}=diag(-1, 1, 1, 1)$, while
for the MCSH system, indices are raised and lowered with respect to the metric $g_{\alpha\beta}=diag(1,-1,-1)$.
We use the convention that the Greek indices such as $\alpha$, $\beta$ run through \{0, 1, 2, 3\} for MKG system,
while they run through \{0, 1, 2\} for MCSH system;
the Latin indices such as $j$, $k$ run through \{1, 2, 3\}, while they run through \{1, 2\} for MCSH system;
and repeated indices are summed.

The corresponding Euler-Lagrange equation of (1) is
\begin{eqnarray}
\partial^{\alpha}F_{\alpha\beta}&=&Im(\phi\overline{D_{\beta}\phi}),\\
D^{\mu}D_{\mu}\phi&=&0.
\end{eqnarray}

Setting $\beta=0$ in the first equation of (3), we obtain the following Gauss-Law constraint
\begin{equation}
\partial_{j}F_{j0}-Im(\phi\overline{D_0\phi})=0.
\end{equation}

The energy of the system (3)-(4) is conserved,
\begin{eqnarray}
E(t)&=&\frac{1}{2}\int_{\mathbb{R}^2}[\sum_{i=1}^{3}F_{0i}^2(x,t)+\sum_{i<j, i, j=1}^{3}F_{ij}^2(x,t)+\sum_{\mu=0}^{3}|D_{\mu}\phi(x,t)|^2]\nonumber\\
&=&E(0), t\geq0.
\end{eqnarray}

The Maxwell-Chern-Simons-Higgs model was proposed in [1] to investigate the self-dual system when there are both
Maxwell and Chern-Simons terms. The corresponding Euler-Lagrange equations are
\begin{eqnarray}
&&\partial_{\lambda}F^{\lambda\rho}+\frac{\kappa}{2}\epsilon^{\mu\nu\rho}F_{\mu\nu}+2e\mathrm{Im}(\phi\overline{D^{\rho}\phi})=0,\nonumber\\
&&D_{\mu}D^{\mu}\phi+U_{\overline{\phi}}(|\phi|^2, N)=0,\\
&&\partial_{\mu}\partial^{\mu}N+U_{N}=0,\nonumber
\end{eqnarray}
where $U(|\phi|^2, N)=\frac{1}{2}(e|\phi|^2+\kappa N-ev^2)^2+e^2N^2|\phi|^2$,
and $U_{\overline{\phi}}, U_{N}$ are formal derivative of $U(|\phi|^2, N)$ with respect to $\overline{\phi}, N$:
\begin{equation}
U_{\overline{\phi}}(|\phi|^2, N)=(e|\phi|^2+\kappa N-ev^2)\phi+e^2N^2\phi,\nonumber
\end{equation}
\begin{equation}
U_{N}(|\phi|^2, N)=\kappa(e|\phi|^2+\kappa N-ev^2)+2e^2N|\phi|^2.\nonumber
\end{equation}

Setting $\rho=0$ in the first equation of (7), we obtain the Gauss-Law constraint
\begin{equation}
\partial_{j}F_{j0}-\kappa F_{12}-2eIm(\phi\overline{D_0\phi})=0.
\end{equation}

The energy of the system (7) is conserved,
\begin{eqnarray}
E(t)&=&\int_{\mathbb{R}^2}[\frac{1}{2}\sum_{i=1}^{2}F_{0i}^2(x,t)+\frac{1}{2}F_{12}^2(x,t)+\sum_{\mu=0}^{2}|D_{\mu}\phi(x,t)|^2\nonumber\\
&&+\sum_{\mu=0}^{2}|\partial_{\mu}N(x,t)|^2+U(|\phi|^2, N)(x,t)]dx=E(0), t\geq0.
\end{eqnarray}

There are two possible boundary conditions to make the energy finite: Either $(\phi, N, A)\rightarrow(0, \frac{ev^2}{\kappa}, 0)$
as $|x|\rightarrow\infty$ or $(|\phi|^2, N, A)\rightarrow (v^2, 0, 0)$ as
$|x|\rightarrow\infty$. The former is called
nontopological boundary condition, and the latter is called topological boundary condition.

For the nontopological boundary condition, we introduce $\widetilde{N}$ satisfying $\widetilde{N}+ev^2/k=N$. Then we have
$(\phi, \widetilde{N}, A_1, A_2)\rightarrow0$ as $|x|\rightarrow\infty$. In this case, $U_{N}$ in the system (7) changes to $U_{\widetilde{N}}$
respectively. For the topological case, we will discuss a subcase of this case,
we assume $\lim_{|x|\rightarrow\infty}\phi=\lambda$ for a fixed complex
scalar $\lambda$ with $|\lambda|=v$, i.e., $\phi$ tends to be constant at the infinity,
this assumption is very natural. We introduce $\varphi$ satisfying $\varphi+\lambda=\phi$. Then we also have $(\varphi, N, A_1, A_2)
\rightarrow 0$ as $|x|\rightarrow\infty$.

The system (3)-(4) are invariant under the gauge transformations
\begin{equation}
A_{\mu}\rightarrow A_{\mu}'=A_{\mu}+\partial_{\mu}\chi,\ \phi\rightarrow \phi'=e^{i\chi}\phi,
\ D_{\mu}\rightarrow D_{\mu}'=\partial_{\mu}-iA_{\mu}'.
\end{equation}
The system (7) is also invariant under the gauge transformations
\begin{equation}
A_{\mu}\rightarrow A_{\mu}'=A_{\mu}+\partial_{\mu}\chi,\ \phi\rightarrow \phi'=e^{ie\chi}\phi,
\ D_{\mu}\rightarrow D_{\mu}'=\partial_{\mu}-ieA_{\mu}'.
\end{equation}
Hence one may impose an additional gauge condition on $A$.
Usually there are three gauge conditions to choose, Coulomb gauge $\partial^i A_i=0$; temporal gauge
$A_0=0$; Lorenz gauge $\partial^{\mu}A_{\mu}=0.$

The Maxwell-Klein-Gordon system is a classical system
which has been studied extensively, see e.g. \cite{A2}, \cite{A3}. For the temporal gauge case, in \cite{A2},
the authors worked in the Coulomb gauge. In this gauge, by exploiting the null structure of the nonlinearity,
they obtained the global existence of finite energy solutions of the system. Then,
by choosing a suitable $\chi$, they use the gauge transform (10) to transform the obtained global solution in the Coulomb gauge to satisfy
the temporal gauge, so they also obtained the global finite energy solution in the temporal gauge.
In this paper, we will work directly on the temporal gauge, and obtain the the global existence of finite energy solutions in this gauge.
We state our results as follows:
\begin{theorem}
Under the temporal gauge $A_0=0$,
given initial data $A_i(0)\in H^{1}(\mathbb{R}^{3}; \mathbb{R})$, $\phi(0)\in H^{1}(\mathbb{R}^{3}; \mathbb{C})$, $\partial_{t}A_i(0)\in L^2(\mathbb{R}^{3}; \mathbb{R})$,
$\partial_{t}\phi(0)\in L^2(\mathbb{R}^{3}; \mathbb{C})$, satisfying the constraint:
\begin{equation}
\partial^{i}(-\partial_{t}A_i(0))=Im(\phi(0)\overline{\partial_t\phi(0)}),\nonumber
\end{equation}
then
\begin{itemize}
\item
(Existence) there exists a global solution $\phi\in C(\mathbb{R}; H^1(\mathbb{R}^3; \mathbb{C}))\cap C^1(\mathbb{R};$ $L^2(\mathbb{R}^3;\\ \mathbb{C}))$, $A_{i}
\in C(\mathbb{R}; H^1(\mathbb{R}^3; \mathbb{R}))\cap C^1(\mathbb{R}; L^2(\mathbb{R}^3; \mathbb{R}))$ satisfying the Maxwell-Klein-Gordon system in the distributional sense.
\item
(Short time uniqueness) If there exists two solutions $\phi_1, \phi_2\in C([-T, T];\\ H^1(\mathbb{R}^3; \mathbb{C}))\cap C^1([-T, T]; L^2(\mathbb{R}^3; \mathbb{C}))$, $A_{1}, A_2
\in C([-T, T]; H^1(\mathbb{R}^3; \mathbb{R}))$$\cap$\\$\ C^1([-T, T]; L^2(\mathbb{R}^3; \mathbb{R}))$ of the system (3)-(4) on the time interval $[-T, T]$ with the same initial data for a sufficiently $T>0$, satisfying for i=1, 2,
\begin{eqnarray}
&&\|(A_i-\nabla(\Delta^{-1}\mathrm{div}A_i(0)))^{\mathrm{df}}\|_{X_{|\tau|=|\xi|}^{\frac{3}{4}+, \frac{3}{4}+}(S_T)}\nonumber\\
&+&\|(A_i-\nabla(\Delta^{-1}\mathrm{divA_i}(0)))^{\mathrm{cf}}\|_{X_{\tau=0}^{1+, \frac{1}{2}+}(S_T)}+\|\phi_i\|_{X_{|\tau|=|\xi|}^{\frac{3}{4}+, \frac{3}{4}+}(S_T)}
<+\infty,\nonumber
\end{eqnarray}
where $S_T:=[-T, T]\times \mathbb{R}^3$, then, $A_1\equiv A_2$, $\phi_1\equiv\phi_2$, a.e. on $S_T$.
\end{itemize}
\end{theorem}
Following the approach in \cite{A4} to investigate the Yang-Mills system in the temporal gauge, we will decompose the
spatial gauge potentials into divergence free parts and curl free parts, and use the $X^{s, b}$ type spaces to
obtain the local well-posedness of the Maxwell-Klein-Gordon system in the temporal gauge.
Here, we see that the estimates for the Maxwell-Klein-Gordon system is similar to that for the Yang-Mills case,
so we can directly use the estimates which have been already proved in \cite{A4}.
To show the finite energy local solution of the system extends globally, as in \cite{A3}, \cite{A5},\footnote{Where, the authors consider the Maxwell-Klein-Gordon system
and Chern-Simons-Higgs system in the Lorenz gauge, while in this gauge, the energy of the corresponding system also can not fully control
the $H^1$ norms of the solutions, and in \cite{A3}, \cite{A5}, the authors use the gauge invariance of these two systems
to transform the solutions such that the initial data satisfy the Coulomb
gauge, and prove the global existence of these two systems
in the finite energy space.} we see that when the initial data satisfies the coulomb gauge, then
their $\dot{H}^1$ norms can be controlled by the energy of the system and $L^2$ norm of the solution,
and also we see in the investigation of the local well-posedness of the system, we have transformed $A$ to $A'$ such that $A'$ satisfy
$(A')^{\mathrm{cf}}=0$, so we have $\mathrm{div} A'=\mathrm{div} (A')^{\mathrm{df}}=0$.
By combining these two facts, we see that the local solution extends globally.

Now we turn our attention to the Cauchy problem of (7). In \cite{A6}, the authors show that
the system is globally well-posed in the Lorenz gauge in $H^{2}\times H^{1}$, and in \cite{A7}, this was extended to $H^{1}\times L^2$ regularity
. Recently, in \cite{A8}, the author investigate the low regularity of the system in the Lorenz gauge.
In \cite{A9}, the authors show that the system is globally well-posed in the temporal gauge in $H^{2}\times H^{1}$. And in \cite{A10}, the present
author show that the system is locally well-posed in the energy regularity $H^1 \times L^2$ and above.

In this paper, by using the approach just described to get the global finite energy solutions of the Maxwell-Klein-Gordon system
in the temporal gauge, we can also get the
global finite energy solutions of the Maxwell-Chern-Simons-Higgs system in the temporal gauge.
Since the nontopological boundary condition is similar to the topological boundary condition case, we just state the results for the
former.
\begin{theorem}
Under the temporal gauge $A_0=0$,
given initial data $A_i(0, x)$, $\widetilde{N}(0, x) \in H^{1}(\mathbb{R}^{2}; \mathbb{R})$, $\phi(0)\in H^{1}(\mathbb{R}^{2}; \mathbb{C})$, $\partial_{t}A_i(0), \partial_t \widetilde{N}(0, x)\in L^2(\mathbb{R}^{2}; \mathbb{R})$, $\partial_{t}\phi(0)\in L^2(\mathbb{R}^{2}; \mathbb{C})$, satisfying the constraint:
\begin{equation}
\partial^{i}(\partial_{t}A_i(0, x))+\kappa F_{12}(0, x)+2eIm(\phi(0)\overline{\partial_t\phi(0)})=0,\nonumber
\end{equation}
then
\begin{itemize}
\item
(Existence) there exists a global solution $\phi\in C(\mathbb{R}; H^1(\mathbb{R}^2; \mathbb{C}))\cap C^1(\mathbb{R}; L^2(\mathbb{R}^2;\\ \mathbb{C}))$, $A_{i}, \widetilde{N}
\in C(\mathbb{R}; H^1(\mathbb{R}^2; \mathbb{R}))\cap C^1(\mathbb{R}; L^2(\mathbb{R}^2; \mathbb{R}))$ satisfying the Maxwell-Chern-Simons-Higgs system in the distributional sense.
\item
(Short time uniqueness) If there exists two solutions $\phi_1, \phi_2 \in C([-T, T];\\ H^1(\mathbb{R}^2;$ $\mathbb{C}))\cap C^1([-T, T]; L^2(\mathbb{R}^2; \mathbb{C}))$, $A_{1}, A_2, \widetilde{N}_1, \widetilde{N}_2
\in C([-T, T]; H^1(\mathbb{R}^2; \mathbb{R})\\ \cap C^1([-T, T]; L^2(\mathbb{R}^2; \mathbb{R}))$
of the system (8) in the temporal gauge $A_0=0$ on the time interval $[-T, T]$ with the same initial data for a sufficiently small $T>0$,
satisfying for $i=1, 2$,
\begin{eqnarray}
&&\|(A_i-\nabla(\Delta^{-1}\mathrm{div}A_i(0)))_{, \pm}^{\mathrm{df}}\|_{X_{|\tau|=|\xi|}^{1, b}(S_T)}\nonumber\\
&+&\|(A_i-\nabla(\Delta^{-1}\mathrm{div}A_i(0)))^{\mathrm{cf}}\|_{X_{\tau=0}^{1+\alpha, b}(S_T)}+\|\phi_{i, \pm}\|_{X_{|\tau|=|\xi|, \pm}^{1, b}(S_T)}\nonumber\\
&+&\|\widetilde{N}_{i, \pm}\|_{X_{|\tau|=|\xi|, \pm}^{1, b}(S_T)}
<+\infty,\nonumber
\end{eqnarray}
for some $b>\frac{1}{2}$, and some sufficiently small $\alpha>0$ and $\delta>0$ such that $\alpha\leq 1-b-\delta$,
where $S_T=[-T, T]\times \mathbb{R}^2$,
then, $A_1\equiv A_2$, $\phi_1\equiv\phi_2$ and $\widetilde{N}_1\equiv \widetilde{N}_2$, a.e. on $S_T$.
\end{itemize}
\end{theorem}

Some notations: Sometimes in the paper we will abbreviate Maxwell-Klein-Gordon as MKG and Maxwell-Chern-Simons-Higgs as MCSH. $H^s(s\in \mathbb{R})$ are Sobolev spaces with respect to the norms $\|f\|_{H^s}=\|\langle\xi\rangle^s \hat{f}\|_{L^2}$, where $\hat{f}
(\xi)=Ff(\xi)$ is the Fourier transform of $f(x)$ and we use the shorthand $\langle\xi\rangle=(1+|\xi|^2)^{\frac{1}{2}}$. We use the shorthand
$X\lesssim Y$ for $X\leq CY$, where $C>>1$ is a constant which may depend on the quantities which are considered fixed. $X\sim Y$ means
$X\lesssim Y\lesssim X$. We use $b+$ to denote $b+\epsilon$, for a sufficiently small positive $\epsilon$, and $\Box:=\partial_{tt}-\Delta$.

In section 2, we will consider the local well-posedness of MKG system in section 2.1, and then
in section 2.2, we will show the global existence of finite energy solution of the system.
In Section 3, we will show the global existence of finite energy solution of the MCSH system.

\section{Maxwell-Klein-Gordon system}
\subsection{Local well-posedness of the MKG system}
Under the temporal gauge $A_0=0$, the system (2)-(3) becomes
\begin{eqnarray}
&&\partial_t\partial^iA_i=-Im(\phi\overline{\partial_t\phi}),\\
&&\Box A_i-\partial_i(\partial^jA_j)=Im(\phi\overline{\partial_i\phi}+iA_i|\phi|^2),\\
&&\Box\phi-i(\partial^iA_i)\phi-2iA^{i}\partial_i\phi-A^iA_i\phi=0.
\end{eqnarray}
Now we decompose $A$ into the divergence free parts $A^{\mathrm{df}}$ and curl free parts $A^{\mathrm{cf}}$. Recall that for a
vector field functions $\vec{X}(x): \mathbb{R}^3\rightarrow \mathbb{R}^3$, $\vec{X}=(-\Delta)^{-1}\mathrm{curlcurl} \vec{X}-(-\Delta)^{-1}\nabla \mathrm{div}\vec{X}$,
Since $\mathrm{div} \mathrm{curl}=0$ and $\mathrm{curl} \nabla=0$, this expresses $\vec{X}$ as the sum of its divergence free and curl free
parts. Let $\mathcal{P}$ denote the projection operator onto the divergence-free vector fields on $\mathbb{R}^3$, $\mathcal{P}:=(-\Delta)^{-1}\mathrm{curl} \mathrm{curl}$.
Then the system (12)-(14) becomes
\begin{eqnarray}
&&\partial_tA^{\mathrm{cf}}=-(-\Delta)^{-1}\nabla[Im(\phi\overline{\partial_t\phi})],\\
&&\Box A^{\mathrm{df}}=-\mathcal{P}[Im(\phi\overline{\partial_i\phi})+iA_i|\phi|^2)],\\
&&\Box\phi-i(\partial^iA_i^{\mathrm{cf}})\phi-2iA_{i}^{\mathrm{df}}\partial^i\phi-2iA_{i}^{\mathrm{cf}}\partial^i\phi-A^iA_i\phi=0.
\end{eqnarray}
We do not expand out $A_i$ in (16)-(17), but remember that $A_i=A_i^{\mathrm{df}}+A_i^{\mathrm{cf}}$.
In [2], the authors show that
\begin{equation}
2A_i^{\mathrm{df}}\partial_i\phi=2A^{\mathrm{df}}\cdot\nabla\phi=Q_{ij}(\phi, |D|^{-1}[R^{i}A^{j}-R^{j}A^{i}]),
\end{equation}
\begin{equation}
\mathcal{P}(\phi\overline{\nabla\phi})_{i}=2R^{j}|D|^{-1}Q_{ij}(Re\phi, Im \phi),
\end{equation}
here $Q_{ij}(u, v):=\partial_iu\partial_jv-\partial_ju\partial_iv, 1\leq i, j\leq3$ denote the null forms,
and $R_i$ are the Riesz transformations defined by $R_i:=|D|^{-1}\partial_i$.
We can also assume that $A^{\mathrm{cf}}(0)=0$. This can be established by using the transform (10), and let $\chi=-\Delta^{-1}\mathrm{div} A(0)$.

Here we construct our solutions in the $X^{s,b}$ type spaces. We recall some definitions and some basic
properties of $X^{s,b}$ and $H^{s,b}$ spaces.

\begin{definition}
For $s, b\in \mathbb{R}$, let $X_{|\tau|=|\xi|, \pm}^{s,b}$ be the completion of the Schwarz space $S(\mathbb{R}^{1+n})$ with respect to
the norm
\begin{equation}
\|u\|_{X_{|\tau|=|\xi|, \pm}^{s,b}}=\|\left<\xi\right>^s\left<-\tau\pm|\xi|\right>^{b}\hat{u}(\tau,\xi)\|_{L_{\tau,\xi}^2},\nonumber
\end{equation}
where $\hat{u}(\tau,\xi)$ denotes the space-time Fourier transformation of $u(t,x)$. Let $X^{s,b}_{|\tau|=|\xi|}$ be the
completion of the Schwarz space $S(\mathbb{R}^{1+n})$ with respect to the norm
\begin{equation}
\|u\|_{X^{s,b}_{|\tau|=|\xi|}}=\|\left<\xi\right>^s\left<|\tau|-|\xi|\right>^{b}\hat{u}(\tau,\xi)\|_{L_{\tau,\xi}^2},\nonumber
\end{equation}
\end{definition}
Clearly, we have
\begin{equation}
\|u\|_{X^{s,b}_{|\tau|=|\xi|}}\leq\|u\|_{X_{|\tau|=|\xi|, \pm}^{s,b}},\ \  \mbox{for}\ \  b\geq0,
\end{equation}
\begin{equation}
\|u\|_{X^{s,b}_{|\tau|=|\xi|}}\geq\|u\|_{X_{|\tau|=|\xi|, \pm}^{s,b}},\ \  \mbox{for}\ \  b\leq0.
\end{equation}
For $S_T=(0, T)\times \mathbb{R}^n$, the restriction space $X_{|\tau|=|\xi|, \pm}^{s,b}(S_T)$ is a Banach space with respect
to the norm
\begin{equation}
\|u\|_{X_{|\tau|=|\xi|, \pm}^{s,b}(S_T)}=\inf\{\|v_{X_{|\tau|=|\xi|, \pm}^{s,b}}:
v\in X_{|\tau|=|\xi|, \pm}^{s,b}\ \ \mbox{ and} \ \ v=u \ \ \mbox{on}\ S_T\}.
\end{equation}
The restriction space $X^{s, b}_{|\tau|=|\xi|}(S_T)$ is defined analogously.
Now we consider the following linear Cauchy problem
\begin{equation}
(-i\partial_{t}\pm\left<\nabla\right>)u=F,\ u|_{t=0}=u_0.
\end{equation}
\begin{lemma}
Let $1/2<b\leq1$, $s\in \mathbb{R}$, $0<T\leq1$, Also, let $0\leq\delta\leq1-b$, Then for $F\in X_{|\tau|=|\xi|, \pm}^{s,b-1+\delta}(S_T)$,
$u_0\in H^s$, the Cauchy problem has a unique solution $u\in X_{|\tau|=|\xi|, \pm}^{s,b}(S_T)$, satisfying the first equation in
the sense of $D'(S_T)$. Moreover,
\begin{equation}
\|u\|_{X_{|\tau|=|\xi|, \pm}^{s,b}(S_T)}\leq C(\|u_0\|_{H^s}+T^{\delta}\|F\|_{X_{|\tau|=|\xi|, \pm}^{s,b-1+\delta}(S_T)}),\nonumber
\end{equation}
where $C$ only depends on $b$.
\end{lemma}
For the linear Cauchy problem
\begin{equation}
\partial_{tt}u-\Delta u=F(t,x), u(0,x)=f(x), \partial_tu(0,x)=g(x),
\end{equation}
we have the following lemma
\begin{lemma}
Let $1/2<b\leq1$, $s\in \mathbb{R}$, $0<T \leq 1$, Also, let $0\leq\delta\leq1-b$. Then for $F\in X^{s-1, b-1+\delta}(S_T)$,
$f\in H^s$, and $g\in H^{s-1}$, there exists a unique $u\in X^{s,b}(S_T)$ solving (2.13) on $S_T$. Moreover,
\begin{equation}
\|u\|_{X^{s,b}(S_T)}\leq C(\|f\|_{H^s}+\|g\|_{H^{s-1}}+T^{\delta}\|F\|_{X^{s-1, b-1+\delta}(S_T)}),
\end{equation}
where $C$ only depends on $b$.
\end{lemma}
\begin{definition}
For $s, b\in \mathbb{R}$, let $X_{\tau=0}^{s,b}$ be the completion of the Schwarz space $S(R^{1+n})$ with respect to
the norm
\begin{equation}
\|u\|_{X_{\tau=0}^{s,b}}=\|\left<\xi\right>^s\left<\tau\right>^{b}\hat{u}(\tau,\xi)\|_{L_{\tau,\xi}^2},\nonumber
\end{equation}
where $\hat{u}(\tau,\xi)$ denotes the space-time Fourier transformation of $u(t,x)$.

For $S_T=(0, T)\times \mathbb{R}^n$, the restriction space $\|u\|_{X_{\tau=0}^{s,b}(S_T)}$ is defined similar to $X_{|\tau|=|\xi|, \pm}^{s,b}(S_T)$
and $X^{s, b}_{|\tau|=|\xi|}(S_T)$.
\end{definition}
Now we consider the following linear Cauchy problem
\begin{equation}
\partial_{t}u=F,\ u|_{t=0}=u_0.
\end{equation}
\begin{lemma}
Let $1/2<b\leq1$, $s\in \mathbb{R}$, $0<T\leq1$, Also, let $0\leq\delta\leq1-b$, Then for $F\in X_{\tau=0}^{s,b-1+\delta}(S_T)$,
$u_0\in H^s$, the Cauchy problem has a unique solution $u\in X_{\tau=0}^{s,b}(S_T)$, satisfying the first equation in
the sense of $D'(S_T)$. Moreover,
\begin{equation}
\|u\|_{X_{\tau=0}^{s,b}(S_T)}\leq C(\|u_0\|_{H^s}+T^{\delta}\|F\|_{X_{\tau=0}^{s,b-1+\delta}(S_T)}),\nonumber
\end{equation}
where $C$ only depends on $b$.
\end{lemma}
We have the following local existence results:
\begin{lemma}
Let $s>\frac{3}{4}$. Given initial data $A_i(0)\in H^{s}(\mathbb{R}^{3}; \mathbb{R})$, $\phi(0)\in H^{s}(\mathbb{R}^{3}; \mathbb{C})$, $\partial_{t}A_i(0)\in H^{s-1}(\mathbb{R}^{3}; \mathbb{R})$,
$\partial_{t}\phi(0)\in H^{s-1}(\mathbb{R}^{3}; \mathbb{C})$, satisfying the constraint:
\begin{equation}
\partial^{i}(-\partial_{t}A_i(0))=Im(\phi(0)\overline{\partial_t\phi(0)}),\nonumber
\end{equation}
then, there exists a time $T>0$, which is a decreasing and continuous function of the data norm
\begin{equation}
\sum_{i=1}^{3}(\|A_i^{\mathrm{df}}(0)\|_{H^s}+\|\partial_tA_i^{\mathrm{df}}(0)\|_{H^{s-1}})+\|\phi(0)\|_{H^s}+\|\partial_t\phi(0)\|_{H^{s-1}},
\end{equation}
and a solution $(A, \phi)$ of (15)-(17) on $(-T, T)\times \mathbb{R}^3$ with the regularity
$$\|A_i^{\mathrm{df}}\|_{X_{|\tau|=|\xi|}^{s, \frac{3}{4}+}(S_T)}+\|A_i^{\mathrm{cf}}\|_{X_{\tau=0}^{s+\frac{1}{4}, \frac{1}{2}+}(S_T)}
+\|\phi\|_{X_{|\tau|=|\xi|}^{s, \frac{3}{4}+}(S_T)}<+\infty.$$
\end{lemma}
We will take the contraction spaces as
\begin{eqnarray}
\|A\|_{X}:&=&\|A^{\mathrm{df}}\|_{X_{|\tau|=|\xi|}^{s, \frac{3}{4}+}(S_T)}+\|A^{\mathrm{cf}}\|_{X_{\tau=0}^{s+\frac{1}{4}, \frac{1}{2}+}(S_T)},\nonumber\\
\|\phi\|_{X}:&=&\|\phi\|_{X_{|\tau|=|\xi|}^{s, \frac{3}{4}+}(S_T)}.\nonumber
\end{eqnarray}
And by using the Lemma 2.3 and Lemma 2.5,
it is standard that the proof of Lemma 2.6 is reduced to the following estimates.
\begin{eqnarray}
\|\nabla^{-1}(\phi\partial_t\psi)\|_{X_{\tau=0}^{s+\frac{1}{4}, -\frac{1}{2}+}}&\lesssim&\|\phi\|_{X_{|\tau|=|\xi|}^{s, \frac{3}{4}+}}
\|\psi\|_{X_{|\tau|=|\xi|}^{s, \frac{3}{4}+}}.
\end{eqnarray}
For $1\leq i, j\leq3$,
\begin{eqnarray}
\|Q_{ij}(|D|^{-1}\phi, \psi)\|_{X_{|\tau|=|\xi|}^{s-1, -\frac{1}{4}+}}&\lesssim&\|\phi\|_{X_{|\tau|=|\xi|}^{s, \frac{3}{4}+}}
\|\psi\|_{X_{|\tau|=|\xi|}^{s, \frac{3}{4}+}},\\
\||D|^{-1}Q_{ij}(\phi, \psi)\|_{X_{|\tau|=|\xi|}^{s-1, -\frac{1}{4}+}}&\lesssim&\|\phi\|_{X_{|\tau|=|\xi|}^{s, \frac{3}{4}+}}
\|\psi\|_{X_{|\tau|=|\xi|}^{s, \frac{3}{4}+}}.\\
\|A_1\phi\psi\|_{X_{|\tau|=|\xi|}^{s-1, -\frac{1}{4}+}}&\lesssim&\|A_1\|_{X_{|\tau|=|\xi|}^{s, \frac{3}{4}+}}\|\phi\|_{X_{|\tau|=|\xi|}^{s, \frac{3}{4}+}}\|\psi\|_{X_{|\tau|=|\xi|}^{s, \frac{3}{4}+}},\\
\|A_1\phi\psi\|_{X_{|\tau|=|\xi|}^{s-1, -\frac{1}{4}+}}&\lesssim&\|A_1\|_{X_{\tau=0}^{s+\frac{1}{4}, \frac{1}{2}+}}\|\phi\|_{X_{|\tau|=|\xi|}^{s, \frac{3}{4}+}}\|\psi\|_{X_{|\tau|=|\xi|}^{s, \frac{3}{4}+}},\\
\|\nabla A^{\mathrm{cf}}\phi\|_{X_{|\tau|=|\xi|}^{s-1, -\frac{1}{4}+}}&\lesssim&\|A^{\mathrm{cf}}\|_{X_{\tau=0}^{s+\frac{1}{4}, \frac{1}{2}+}}
\|\phi\|_{X_{|\tau|=|\xi|}^{s, \frac{3}{4}+}},\\
\|A^{\mathrm{cf}}\nabla\phi\|_{X_{|\tau|=|\xi|}^{s-1, -\frac{1}{4}+}}&\lesssim&\|A^{\mathrm{cf}}\|_{X_{\tau=0}^{s+\frac{1}{4}, \frac{1}{2}+}}
\|\phi\|_{X_{|\tau|=|\xi|}^{s, \frac{3}{4}+}},
\end{eqnarray}
\begin{eqnarray}
\|A_1A_2\phi\|_{X_{|\tau|=|\xi|}^{s-1, -\frac{1}{4}+}}&\lesssim&\prod_{i=1}^2\min(\|A_i\|_{X_{|\tau|=|\xi|}^{s, \frac{3}{4}+}},
\|A_i\|_{X_{\tau=0}^{s+\frac{1}{4}, \frac{1}{2}+}})\|\phi\|_{X_{|\tau|=|\xi|}^{s, \frac{3}{4}+}}.\nonumber\\
\end{eqnarray}
These estimates have appeared in [4], please see the estimates (15)-(19) in \cite{A4}, we refer the reader to the proof of these
estimates in this paper, and we complete the proof of Lemma 2.6.
\subsection{Global existence of MKG system}
Now we turn to the global existence part of Theorem 1.1. We will work on the finite energy space of $(A, \phi)$,
i.e. we will work on $H^1\times L^2$ regularity for them.
Suppose $(A, \phi)$ are the solutions of the system (12)-(14) on $[0, T]$, then we have on $[0, T]$,
\begin{equation}
\frac{d}{dt}\frac{1}{2}\int|A(t)|^2dx=\int A\cdot\partial_tAdx\leq\|A(t)\|_{L^2}\|\partial_tA(t)\|_{L^2},
\end{equation}
thus $\frac{d}{dt}\|A(t)\|_{L^2}\leq\|\partial_tA(t)\|_{L^2}$, whence
\begin{equation}
\|A(t)\|_{L^2}\leq\|A(0)\|_{L^2}+\int_0^t\|\partial_tA(t')\|_{L^2}dt'\leq\|A(0)\|_{L^2}+Ct\sqrt{E(0)}.
\end{equation}
Also, we have
\begin{equation}
\frac{d}{dt}\frac{1}{2}\int|\phi(t)|^2dx=Re\int\phi\overline{\partial_t\phi}dx=Re\int\phi\overline{D_0\phi}dx\leq\|\phi\|_{L^2}\sqrt{E(0)},
\end{equation}
thus we have
\begin{equation}
\|\phi(t)\|_{L^2}\leq\|\phi(0)\|_{L^2}+Ct\sqrt{E(0)}.
\end{equation}
By combining (37) and (39), we have
\begin{lemma}
Suppose $(A_i, \phi)\in C([-T, T]; H^1(\mathbb{R}^3; \mathbb{R}))\cap C^{1}([-T, T]; L^2(\mathbb{R}^3;\\ \mathbb{R}))$$\times C([-T, T]; H^1(\mathbb{R}^3; \mathbb{C}))\cap C^{1}([-T, T]; L^2(\mathbb{R}^3; \mathbb{C}))$ be the solutions of the MKG system
(12)-(14) on the time interval $[-T, T]$ under the temporal
gauge $A_0=0$, then $\forall\ t\in[-T, T], \|A(t)\|_{L^2}+\|\phi(t)\|_{L^2}\leq C(E(0), T)$, where $C(E(0), T)$ denotes a function depending on $E(0)$ and $T$.
\end{lemma}
Now we turn to estimate the $\dot{H}^1$ norms of $A$ and $\phi$, and we recall a Lemma from [3].
\begin{lemma}
Suppose we are given $\phi_0\in L^2(\mathbb{R}^3; \mathbb{C})$ and $\vec{a}\in\dot{H}^1(\mathbb{R}^3; \mathbb{R}^3)$. Define
\begin{equation}
U=\nabla\phi_0-i\phi_0\vec{a},
\end{equation}
and assume that $U\in L^2$. Then $\phi_0\in\dot{H}^1$, and
\begin{equation}
\|\nabla\phi_0\|_{L^2}\leq2\|U\|_{L^2}+C\|\vec{a}\|_{\dot{H}^1}^2\|\phi_0\|_{L^2},
\end{equation}
where $C$ is an absolute constant.
\end{lemma}
Also for vector field functions $\vec{a}(x)\in \mathbb{R}^3\rightarrow \mathbb{R}$, if $\nabla\times\vec{a}=B$ and $\nabla\cdot\overrightarrow{a}
=0$, then $\|\vec{a}\|_{\dot{H}^1}\leq C\|B\|_{L^2}$. Since $\partial^iA_i^{\mathrm{df}}=0$, we have
\begin{equation}
\|A^{\mathrm{df}}\|_{\dot{H}^1}\leq C\|\mathrm{curl} A^{\mathrm{df}}\|_{L^2}.
\end{equation}
So, by using the Lemma 2.8 for $U=D_i\phi$ for $i=1, 2, 3$ and (42), when $A^{\mathrm{cf}}(0)=0$, we have
\begin{eqnarray}
\|\nabla\phi(0)\|_{L^2}&\leq&C\sum_{i=1}^{3}\|D_i\phi(0)\|_{L^2}+C\|A^{\mathrm{df}}(0)\|_{\dot{H}^1}^2\|\phi(0)\|_{L^2}\nonumber\\
&\leq&C(1+E(0)(1+\|\phi(0)\|_{L^2})).
\end{eqnarray}
Under the temporal gauge $A_0=0$, we have $F_{0i}=\partial_{t}A_i$, $D_0\phi=\partial_t\phi$, so
\begin{equation}
\|\partial_{t}A_i\|_{L^2}\leq\sqrt{E(t)},\ \|\partial_t\phi\|_{L^2}\leq\sqrt{E(t)}.
\end{equation}

Now for the initial data $(A_i(0), \phi(0))$ of the system (12)-(14), firstly we make a gauge transform of the form (6)
with $\chi=-\Delta^{-1}\mathrm{div} A(0)$ to let the transformed initial data $A'$ satisfy $(A')^{\mathrm{cf}}(0)=0$.
For the transformed initial data $(A'(0), \phi'(0))$, since $A'(0)=(A'(0))^{\mathrm{df}}+(A'(0))^{\mathrm{cf}}=(A'(0))^{\mathrm{df}}$,
and note that the energy $E(t)$ is invariant under the transform (10),
so by Lemma 2.7, (42)-(44), the corresponding ${\dot{H}}^1$ norm of $(A'(0), \phi'(0))$, $\|\partial_{t}A_i\|_{L^2}$, and
$\|\partial_t\phi\|_{L^2}$ can be controlled by the $L^2$ norm of $(A(0), \phi(0))$ and $E(0)$.
Also for the $L^2$ norms of $A'(0)$ and $\phi'(0)$, we have
\begin{equation}
\|\phi'(0)\|_{L^2}=\|\phi(0)\|_{L^2},
\end{equation}
and
\begin{eqnarray}
&&\|A'(0)\|_{L^2}\nonumber\\
&=&\|A(0)+\partial_i\chi\|_{L^2}=\|A(0)-\partial_i\Delta^{-1}\mathrm{div} A(0)\|_{L^2}\nonumber\\
&\leq&\|A(0)\|_{L^2}+\|\partial_i\Delta^{-1}\mathrm{div} A(0)\|_{L^2}\leq C\|A(0)\|_{L^2}.
\end{eqnarray}
Now use the
Lemma 2.6 for the transformed initial data $(A'(0), \phi'(0))$, we can get a local solution on $[0, T_1]$ with $T_1$ depending on the
$\|A(0)\|_{L^2}$,  $\|\phi(0)\|_{L^2}$ and $E(0)$,
then we use the inverse transformation
of the form (6) with $\chi=\Delta^{-1}\mathrm{div} A(0)$
to transform back to get the solution for the original initial data, and thus we obtain the solution of the system on $[0, T_1]$. Then, we
begin with $(A(T_1), \phi(T_1))$ as the initial data for the system (12)-(14), for the transformed initial data $(A'(T_1), \phi'(T_1))$, since $A'(T_1)=(A'(T_1))^{\mathrm{df}}+(A'(T_1))^{\mathrm{cf}}=(A'(T_1))^{\mathrm{df}}$,
similarly to get the solution on $[0, T_1]$, by using Lemma 2.7, (42)-(44), and modifying (45)-(46) with $(\phi'(0), A'(0))$
replaced by $(\phi'(T_1), A'(T_1))$, we obtain
a solution of the system on $[T_1, T_2]$ with the length of $[T_1, T_2]$ depending on $\|A(T_1)\|_{L^2}$,  $\|\phi(T_1)\|_{L^2}$ and $E(0)$.
We can repeat this procedure to obtain the solution of the system on $[0, T^{*})=\bigcup_{i\geq1}
[T_{i-1}, T_i]$, with $T_0=0$ and the length of $[T_{i-1}, T_i]$ depending on $\|A(T_{i-1})\|_{L^2}$,  $\|\phi(T_{i-1})\|_{L^2}$ and $E(0)$.

Now we need to show $T^{*}=+\infty$. For this, we show that we can obtain the solution on $[0, T_{*}]$ for any fixed $0<T_{*}<+\infty$. By
Lemma 2.7, we see that for $i\geq1$, when $T_{i-1}\leq T_{*}$, $\|A(T_{i-1})\|_{L^2}$,  $\|\phi(T_{i-1})\|_{L^2}$ can be bounded by a function depends on $\|A(0)\|_{L^2}$,  $\|\phi(0)\|_{L^2}$, $E(0)$, $T_{*}$, thus on $[0, T_{*}]$, the length of each existence interval
in the iteration process only depends on $\|A(0)\|_{L^2}$,  $\|\phi(0)\|_{L^2}$, $E(0)$, $T_{*}$, and can
be chosen independently of the iteration step,
so after finitely
many iteration steps one obtains a solution on any interval $[0, T_{*}]$, and prove $T^{*}=+\infty$,
so we get a global solution $(A, \phi)$ of the Maxwell-Klein-Gordon system on $[0, +\infty)$. Similarly, we
can get a global solution $(A, \phi)$ of the Maxwell-Klein-Gordon system on $(-\infty, 0]$. By combining these, we get
a global solution $(A, \phi)$ of the Maxwell-Klein-Gordon system on the time interval $(-\infty, +\infty)$ in the finite energy space.
The uniqueness part of Theorem 1.1 is obvious, and we complete the proof of Theorem 1.1.

\section{Maxwell-Chern-Simons-Higgs system}
Under the temporal gauge and nontopological boundary condition, we rewrite the Euler-Lagrange equations (7) in terms of $(A, \phi, \widetilde{N})$ as follows
\begin{eqnarray}
\partial_t(\mathrm{div} A)&+&\kappa(\partial_1A_2-\partial_2A_1)+2e\mathrm{Im}(\phi\partial_t\bar{\phi})=0,\nonumber\\
\Box A_1&+&\partial_1(\mathrm{div} A)+\kappa\partial_tA_2+2e\mathrm{Im}(\phi\partial_1\bar{\phi})+2e^2A_1|\phi|^2=0,\nonumber\\
\Box A_2&+&\partial_2(\mathrm{div} A)-\kappa\partial_tA_1+2e\mathrm{Im}(\phi\partial_2\bar{\phi})+2e^2A_2|\phi|^2=0,\nonumber\\
\Box \phi&=&-2ieA_{j}\partial_{j}\phi-ie\partial_jA_j\phi-e^2A_{j}^2\phi-(e|\phi|^2+\kappa N-ev^2)\phi-e^2N^2\phi,\nonumber\\
\Box \widetilde{N}&=&-\kappa(e|\phi|^2+\kappa N-ev^2)-2e^2N|\phi|^2,
\end{eqnarray}
with given initial data $A_{i}(0, x), \phi(0,x), \widetilde{N}(0,x),\partial_{t}A_{i}(0,x)$,$\partial_{t}\phi(0,x)$,
$\partial_{t}\widetilde{N}(0,x)$ satisfying the constraint
\begin{equation}
\partial_{j}\partial_{t}A_{j}(0, x)+\kappa F_{12}(0, x)+2e\mathrm{Im}(\phi\partial_t\bar{\phi})(0, x)=0.
\end{equation}
We do not expand out $N$ in the right hand side of system (47), but remember that $N=\widetilde{N}+\frac{ev^2}{k}$.

We decompose $A=(A_1, A_2)$ into divergence free part $A^{\mathrm{df}}=(A^{\mathrm{df}}_1, A^{\mathrm{df}}_2)$ and curl free part
$A^{\mathrm{cf}}=(A^{\mathrm{cf}}_1, A^{\mathrm{cf}}_2)$, such that
$A=A^{\mathrm{df}}+A^{\mathrm{cf}}$, $\mathrm{div} A^{\mathrm{df}}=0$, and $\mathrm{curl} A^{\mathrm{cf}}=0$. Remember that for $u=(u_1(x, y), u_2(x, y))$, $\mathrm{div} u(x,y)=u_{1,x}(x, y)+u_{2, y}(x, y)$, $\mathrm{curl} u(x, y)=u_{2, x}(x, y)-u_{1, y}(x, y)$. And by calculating,
we have $A^{\mathrm{df}}=(-\partial_2\Delta^{-1}(\frac{\partial A_2}{\partial x}-\frac{\partial A_1}{\partial y}),
\partial_1\Delta^{-1}(\frac{\partial A_2}{\partial x}-\frac{\partial A_1}{\partial y}))$, $A^{\mathrm{cf}}=\Delta^{-1}(\nabla \mathrm{div} A).$
\begin{eqnarray}
\partial_t A^{\mathrm{cf}}&=&-\Delta^{-1}\nabla[\kappa(\partial_1A_2-\partial_2A_1)+2e\mathrm{Im}(\phi\partial_t\bar{\phi})],\nonumber\\
\Box A_1^{\mathrm{df}}&=&-\kappa\Delta^{-1}\partial_{12}\partial_tA_1-\kappa\Delta^{-1}\partial_{22}\partial_tA_2
+2e\Delta^{-1}\partial_{12}\mathrm{Im}(\phi\partial_2\bar{\phi})\nonumber\\
&&-2e\Delta^{-1}\partial_{22}\mathrm{Im}(\phi\partial_1\bar{\phi})
+2e^2\Delta^{-1}\partial_{12}(A_2|\phi|^2)-2e^2\Delta^{-1}\partial_{22}(A_1|\phi|^2),\nonumber\\
\Box A_2^{\mathrm{df}}&=&\kappa\Delta^{-1}\partial_{11}\partial_tA_1+\kappa\Delta^{-1}\partial_{21}
\partial_tA_2-2e\Delta^{-1}\partial_{11}\mathrm{Im}(\phi\partial_2\bar{\phi})\nonumber\\
&&+2e\Delta^{-1}\partial_{12}\mathrm{Im}(\phi\partial_1\overline{\phi})-2e^2\Delta^{-1}\partial_{11}(A_2|\phi|^2)
+2e^2\Delta^{-1}\partial_{12}(A_1|\phi|^2),\nonumber\\
\Box \phi&=&-2ieA_{j}\partial_{j}\phi-2ie\partial_jA_j\phi-e^2A_{j}^2\phi-(e|\phi|^2+\kappa N-ev^2)\phi-e^2N^2\phi,\nonumber\\
\Box \widetilde{N}&=&-\kappa(e|\phi|^2+\kappa N-ev^2)-2e^2N|\phi|^2.
\end{eqnarray}
We can also assume
\begin{equation}
A^{\mathrm{cf}}(0)=0.
\end{equation}
This can be established by using the transformation (1.11),
and let\\ $\chi=-\Delta^{-1}\mathrm{div}A(0)$.

In [10], by setting
\begin{eqnarray}
&&A_{i}^{\mathrm{df}}=A_{i,+}^{\mathrm{df}}+A_{i,-}^{\mathrm{df}}, \phi=\phi_{+}+\phi_{-}, \widetilde{N}=\widetilde{N}_{+}+\widetilde{N}_{-},\nonumber\\
&&A_{i,\pm}^{\mathrm{df}}=\frac{1}{2}(A_{i}^{\mathrm{df}}\pm i^{-1}\left<\nabla\right>^{-1}\partial_tA_{i}^{\mathrm{df}}), i=1,2,
\phi_{\pm}=\frac{1}{2}(\phi\pm i^{-1}\left<\nabla\right>^{-1}\partial_{t}\phi),\nonumber\\
&&\widetilde{N}_{\pm}=\frac{1}{2}(\widetilde{N}\pm i^{-1}\left<\nabla\right>^{-1}\partial_{t}\widetilde{N}),\nonumber
\end{eqnarray}
and consider the equations satisfied by $A_{i}^{\mathrm{cf}}$, $A_{i,\pm}^{\mathrm{df}}$, $\phi_{\pm}, \widetilde{N}_{\pm}$,
we obtain the local well-posedness result of the MCSH in the temporal gauge, see Theorem 3.2 and Theorem 3.1 in \cite{A10}.
From which, we can deduce the following local well-posedness results for the system (49) under the conditions (48) and (50).
\begin{theorem}
The Maxwell-Chern-Simons-Higgs Cauchy problem (49), (48), (50) is locally well-posed in $H^s\times H^{s-1}$,
$s\geq1$, in nontopological boundary conditions. To be precise, there exists a time $T$
depending on the initial data norms
$\|A_i^{\mathrm{df}}(0)$, $\phi(0)$, $\widetilde{N}(0)\|_{H^s\times H^{s-1}}$, and a solution $(A_i^{\mathrm{cf}}(t), A_i^{\mathrm{df}}(t), \phi(t), \widetilde{N}(t)+ev^2/k)$ of (49), (48) and (50) on $[-T, T]\times \mathbb{R}^2$ with the regularity
\begin{eqnarray}
\phi_{\pm}, \widetilde{N}_{\pm},A_{i,\pm}^{\mathrm{df}}\in X_{_{|\tau|=|\xi|}, \pm}^{s,b}(S_T)\subset C([-T,T],H^s),\nonumber\\
A_{i}^{\mathrm{cf}}\in X_{\tau=0}^{s+\alpha,b}(S_T)\subset C([-T,T],H^{s+\alpha}),
s\geq1,\nonumber
\end{eqnarray}
for some $b>\frac{1}{2}$, and some sufficiently small $\alpha>0$ and $\delta>0$ such that
$\alpha\leq1-b-\delta$, where $S_T=[-T, T]\times \mathbb{R}^2$. The solution is unique in this regularity class.
\end{theorem}
Now we turn to the global existence part of Theorem 1.1. We will work on the finite energy space of $(A, \phi, \widetilde{N})$,
i.e. we will work on $H^1\times L^2$ regularity for them.
Similar to the Maxwell-Klein-Gordon system, we can control the $L^2$ norm of $A(t)$ and $\phi(t)$. For the
$L^2$ norms of $\widetilde{N}(t)$, we see by the expression of the energy $E(t)$, we have
\begin{equation}
\|\partial_t\widetilde{N}(t)\|_{L^2}+\|\nabla \widetilde{N}(t)\|_{L^2}\leq CE(0)^{\frac{1}{2}},
\end{equation}
so we have
\begin{eqnarray}
&&\|\widetilde{N}(t)\|_{L^2}\nonumber\\
&=&\|\widetilde{N}(0)+\int_0^t\|\partial_t\widetilde{N}(t')dt'\|_{L^2}\leq\|\widetilde{N}(0)\|_{L^2}
+\int_0^t\|\partial_t\widetilde{N}(t')\|_{L^2}dt'\nonumber\\
&\leq&\|\widetilde{N}(0)\|_{L^2}+CE(0)^{\frac{1}{2}}t.
\end{eqnarray}
By combining these, we have
\begin{lemma}
Suppose $(A, \phi, \widetilde{N})\in C([-T, T]; H^1(\mathbb{R}^2; \mathbb{R}))\cap C^{1}([-T, T]; L^2(\mathbb{R}^2; \mathbb{R})) \times C([-T, T]; H^1(\mathbb{R}^2; \mathbb{R}))\cap C([-T, T]; L^2(\mathbb{R}^2; \mathbb{R}))\times C([-T, T]; H^1(\mathbb{R}^2; \mathbb{C}))\\ \cap\ C^{1}([-T, T]; L^2(\mathbb{R}^2; \mathbb{C}))$,
and $(A, \phi, \widetilde{N}+ev^2/k)$ be the solutions of the MCSH system under the temporal gauge $A_0=0$, then $\forall\ t\in[-T, T]$, $\|A(t)\|_{L^2}+\|\phi(t)\|_{L^2}+\|\widetilde{N}(t)\|_{L^2}\leq C(E(0), T)$, where $C(E(0), T)$ denotes a function depending on $E(0)$ and $T$.
\end{lemma}
Now for vector field functions $(A_1(x), A_2(x))$, if $\partial_1A_2-\partial_2A_1=B$, $\partial_1A_1+\partial_2A_2=0$, then
we have $\|\nabla A\|_{L^2}\leq\|B\|_{L^2}$. So by using this inequality for $A^{\mathrm{df}}$, we have
\begin{equation}
\|\nabla A^{\mathrm{df}}\|_{L^2}\leq\|\partial_1A_2^{\mathrm{df}}-\partial_2A_1^{\mathrm{df}}\|_{L^2}.
\end{equation}
Also, for $i=1, 2$, since $D_i\phi(0)=\partial_i\phi(0)-i\phi(0)A_i^{\mathrm{df}}$, we have
\begin{eqnarray}
&&\|\partial_i\phi(0)\|_{L^2}\nonumber\\
&\leq&\|D_i\phi(0)\|_{L^2}+\|i\phi(0)A_i^{\mathrm{df}}\|_{L^2}
\leq\|D_i\phi(0)\|_{L^2}+\|\phi(0)\|_{L^4}\|A_i^{\mathrm{df}}\|_{L^4}\nonumber\\
&\leq&\|D_i\phi(0)\|_{L^2}+\|\phi(0)\|_{L^2}^{\frac{1}{2}}\|\nabla\phi(0)\|_{L^2}^{\frac{1}{2}}\|
A_i^{\mathrm{df}}\|_{L^2}^{\frac{1}{2}}\|\nabla A_i^{\mathrm{df}}\|_{L^2}^{\frac{1}{2}}\nonumber\\
&\leq&\|D_i\phi(0)\|_{L^2}+\frac{1}{2}\|\nabla\phi(0)\|_{L^2}+2\|\phi(0)\|_{L^2}\|A^{\mathrm{df}}\|_{L^2}\|\mathrm{curl} A^{\mathrm{df}}\|_{L^2}.
\end{eqnarray}
So,
\begin{equation}
\|\nabla\phi(0)\|_{L^2}\leq C+CE(0)+C\|\phi(0)\|_{L^2}^2\|A^{\mathrm{df}}\|_{L^2}^2.
\end{equation}
By combining (51), (53) and (55),
we see that the $H^1$ norms of $(A, \phi, \widetilde{N})$ can be controlled by $E(0)$ and their $L^2$ norms,
so begin with the initial data, we can proceed as in the Maxwell-Klein-Gordon case to get a global finite energy solution of the Maxwell-Chern-Simons-Higgs system. The uniqueness part of Theorem 1.2 is obvious, and we complete the proof of Theorem 1.2.

\section*{Acknowledgments}  The author would like to deeply
thank the referees for many invaluable comments and suggestions which improve this paper.




\end{document}